\documentclass [10pt] {article}

 \usepackage [english] {babel}
 \usepackage {amssymb}
 \usepackage {amsmath}
 \usepackage {amsthm}
 \usepackage [curve] {xypic}

 \let \tildeaccent=\~

 \DeclareTextCommandDefault {\cprime}{\/{\mathsurround=0pt$'$}}

 \def \~{\tilde}
 \def \:{\colon}
 
 \def \ge{\geqslant}
 \def \<{\langle}
 \def \>{\rangle}
 \def \smash{\wedge}
 \def \o{\circ}
 \def \x{\times}
 \def \d{\partial}
 \def \={\overline}
 \def \_{\underline}
 \def \^#1{{}^{#1}}

 \def \Z{\boldsymbol{\mathsf Z}}
 \def \id{\mathrm{id}}

 \DeclareMathOperator {\Hom}{Hom}
 \DeclareMathOperator {\diag}{diag}
 \DeclareMathOperator {\Tot}{Tot}

 \newcommand* {\Mod} [1]
 {\text{\rm$#1$-Mod}}

 \newcommand* {\head} [1]
 {\subsubsection * {#1}}

 \newcommand* {\subhead} [1]
 {\medskip \noindent {\bf\itshape #1\/}}

 \newenvironment* {claim} [1] []
 {\begin{trivlist}\item [\hskip\labelsep {\bf #1}] \it}
 {\end{trivlist} }

 \newenvironment* {demo} [1] []
 {\begin{trivlist}\item [\hskip\labelsep {\it #1}] }
 {\end{trivlist} }

 \newenvironment* {remark} [1] []
 {\begin{trivlist}\item [\hskip\labelsep {\it #1\/}] }
 {\end{trivlist} }

\begin {document}

 \title {\Large\bf On homology of map spaces}

 \author {\normalsize\rm
          S.~S.~Podkorytov}

 \date {}

 \maketitle

 \begin {abstract} \noindent
 Following an idea of Bendersky--Gitler,
 we construct an isomorphism between Anderson's and Arone's
 complexes modelling the chain complex of a map space.
 This allows us to apply Shipley's convergence theorem to
 Arone's model.
 As a corollary,
 we reduce
 the problem of homotopy equivalence for certain ``toy'' spaces
 to a problem in homological algebra.
 \end {abstract}


 \head {}


 A {\it space\/} is a pointed simplicial set.
 A {\it map\/} is a basepoint-preserving simplicial map.
 Chains, homology etc.\ are reduced
 with coefficients in a commutative ring $R$.

 Fix spaces $X$ and $Y$.
 We are interested in the homology of $Y^X$,
 the space of maps $X\to Y$.


 \subhead {0.A. Arone's approach.}
 Let $\Omega$ be the category
 whose objects are the sets $\<s\>=\{1,\dotsc,s\}$, $s>0$, and
 whose morphisms are surjective functions.
 Let $\Omega^\o$ denote the dual category.
 For $n\in\Z$,
 let us define a functor $M_n(X)\:\Omega^\o\to\Mod R$.
 Set $M_n(X)(\_s)=C_n(X^{\smash s})$,
 where $X^{\smash s}$ is the $s$th smash power.
 For a morphism $h\:\<t\>\to\<s\>$,
 set $M_n(X)(h)=C_n(h^\sharp)\:
 C_n(X^{\smash s})\to C_n(X^{\smash t})$, where
 the map $h^\sharp\:X^{\smash s}\to X^{\smash t}$ is given by
 $h^\sharp(x_1\dotso x_s)=x_{h(1)}\dotso x_{h(t)}$
 for $x_1,\dotsc,x_s\in X_n$, $n\ge0$.
 Here the simplex $x_1\dotso x_s\in(X^{\smash s})_n$ is
 the image of the simplex $(x_1,\dotsc,x_s)\in(X^s)_n$
 under the projection.

 \begin {claim} [0.1. Lemma.]
 The functors $M_n(X)$ are projective objects of the abelian
 category of functors $\Omega^\o\to\Mod R$.
 \end {claim}

 Proof is given in 1.B.

 The boundary operators
 $\d\:C_n(X^{\smash s})\to C_{n-1}(X^{\smash s})$ form a
 functor morphism $\d:M_n(X)\to M_{n-1}(X)$.
 Thus $M_*(X)$ is a chain complex of functors.

 \begin {claim} [0.2. Corollary.]
 If a map $e\:X\to Y$ is a weak equivalence,
 then the induced chain homomorphism $M_*(e)\:M_*(X)\to M_*(Y)$
 is a chain homotopy equivalence.
 \qed
 \end {claim}

 We have
 the (unbounded) chain complex of $R$-modules
 $$
 G_*(X,Y)=\Hom_*(M_*(X),M_*(Y))
 $$
 and
 a chain homomorphism
 $$
 \lambda_*(X,Y)\:C_*(Y^X)\to G_*(X,Y),
 $$
 see 2.C, 2.D.
 A natural filtration of $G_*(X,Y)$ yields
 the Arone spectral sequence
 \begin{equation}
 H_{t-s}(\Hom_{\Sigma_s\,*}(C_*(X^{(s)}),C_*(Y^{\smash s})))=
 \^1E^s_t\Rightarrow
 H_{t-s}(G_*(X,Y)),
 \tag{$*$}
 \end{equation}
 where $X^{(s)}=X^{\smash s}/(\text{\rm fat diagonal})$
 \cite{Aro}, \cite{AheKuh}.
 \cite[Theorem 9.2]{Boa} ensures conditional convergence.
 If $Y$ is $(\dim X)$-connected, then
 the convergence is strong and
 $\lambda_*(X,Y)$ is a quasi-isomorphism,
 see \cite{Aro} for the precise statement.
 (A similar result was obtained in \cite[Ch.\ III, \S~5]{Vas}.)
 We wish to get free of the connectivity assumption.
 

 \subhead {0.B. Main results.}
 Here we suppose $R=\Z/\ell$, $\ell$ a prime.
 We call $Y$ $\ell$-toy if
 $\pi_0(Y)$ is finite and
 $\pi_n(Y,y)$ is a finite $\ell$-group
 for all $y\in Y_0$ and $n>0$.

 \begin {claim} [0.3. Theorem.]
 Suppose that
 $X$ is essentially compact\footnote
 {A space is {\it compact\/} (or {\it finite\/}) if
  it is generated by a finite number of simplices.
  {\it Essentially compact\/} means
  weakly equivalent to a compact space.}
 and
 $Y$ is fibrant and $\ell$-toy.
 Then
 $\lambda_*(X,Y)$ is a quasi-isomorphism.
 \end {claim}

 This follows from Theorems 0.5 and 0.6 below,
 see \S~4 for details.
 Under the assumptions of the theorem,
 the convergence of ($*$) is strong
 by \cite[Theorem 7.1]{Boa}.

 \begin {claim} [0.4. Corollary.]
 Suppose that
 $X$ and $Y$ are essentially compact and $\ell$-toy.
 Suppose that
 the complexes $M_*(X)$ and $M_*(Y)$ are chain homotopy
 equivalent.
 Then
 $X$ and $Y$ are weakly equivalent.
 \end {claim}

 The proof is given in \S~5.
 There seems to be no easy/functorial way to extract $\pi_1(X)$
 or the ring structure of $H^*(X)$ from $M_*(X)$.
 The corollary has an algebraic analogue \cite{Pod}.


 \subhead {0.C. Anderson's approach.}
 For a pointed set $S$,
 the space $Y^S$ is defined to be the fibre of the projection
 $$
 \prod_{s\in S}Y\to Y
 $$
 corresponding to $s=*$
 (this agrees with our convention that maps preserve
 basepoints).

 We have
 an (unbounded) chain complex $D_*(X,Y)$ with
 $$
 D_n(X,Y)=\prod_{q-p=n}C_q(Y^{X_p})
 $$
 and a chain homomorphism
 $$
 \mu_*(X,Y)\:C_*(Y^X)\to D_*(X,Y),
 $$
 see 2.F, 2.G for details.
 A natural filtration of $D_*(X,Y)$ yields
 the Anderson spectral sequence
 $$
 H_q(Y^{X_p})=
 \^1E^p_q\Rightarrow
 H_{q-p}(D_*(X,Y)).
 $$
 If $Y$ is $(\dim X)$-connected, then
 $\mu_*(X,Y)$ is a quasi-isomorphism,
 see \cite{And} and \cite[4.2]{Bou} for precise statements.
 Shipley got rid of the connectivity assumption \cite{Shi}.

 \begin {claim} [0.5. Theorem.]
 Suppose that $R=\Z/\ell$, $\ell$ a prime.
 Suppose that
 $X$ is compact and
 $Y$ is fibrant and $\ell$-toy.
 Then
 $\mu_*(X,Y)$ is a quasi-isomorphism.
 \end {claim}

 This is a special case of Shipley's strong convergence
 theorem,
 see \S~3 for details.


 \subhead {0.D. Comparing $G_*(X,Y)$ and $D_*(X,Y)$.}
 We construct a chain homomorphism
 $$
 \epsilon_*(X,Y)\:D_*(X,Y)\to G_*(X,Y)
 $$
 such that the diagram
 $$
 \xymatrix {
 & &
 D_*(X,Y)
 \ar[dd]^-{\epsilon_*(X,Y)} \\
 C_*(Y^X)
 \ar[rru]^-{\mu_*(X,Y)}
 \ar[rrd]_-{\lambda_*(X,Y)} & & \\
 & &
 G_*(X,Y)
 }
 $$
 is commutative,
 see 2.H.

 \begin {claim} [0.6. Theorem.]
 Suppose that $X$ is gradual\footnote
 {A space $X$ is {\it gradual\/} (or {\it finite type\/}) if
 the sets $X_n$, $n\ge0$, are finite.}.
 Then $\epsilon_*(X,Y)$ is an isomorphism.
 \end {claim}

 Proof is given in 2.I.

 \begin {remark} [Remark.]
 In some cases,
 the ${}^2E$ term of the Anderson spectral sequence
 \cite[Theorem 7.1 (2)]{BenGit} and
 the ${}^1E$ term of the Arone spectral sequence
 differ in the grading only.
 This suggested relation of the two approaches
 \cite[footnote 1]{AheKuh}
 and motivated this work.
 Our construction of $\epsilon_*(X,Y)$ follows
 the line of \cite[\S~6]{BenGit}.
 \end {remark}


 \subhead {Acknoledgement.}
 I am grateful to
 S.~Betley,
 V.~A.~Vassiliev and
 M.~Yu.~Zva\-gel{\cprime}\-s\-ki\u{\i}
 for useful discussions.


 \head {1. Preliminaries}


 \subhead {1.A. Notation.}
 For a pointed set $S$,
 we put $S^\x=S\setminus\{*\}$.

 $\Delta^p_+$ is the standard $p$-simplex with an added
 basepoint.
 Let $\iota_p\in(\Delta^p_+)_p$ be the fundamental simplex.

 For $x\in X_n$,
 $[x]\in C_n(X)$ is the chain consisting of
 the single simplex $x$ with the coefficient $1$.

 Given functors $F,F'\:\Omega^\o\to\Mod R$,
 a functor morphism $T\:F\to F'$ consists of homomorphisms
 $\^sT\:F(\<s\>)\to F'(\<s\>)$.

 \subhead {1.B. Proof of Lemma 0.1}
 (cf. \cite[\S~I]{AntBet}).
 Fix a linear order on $X_n^\x$.
 Introduce the set
 $$
 I=\coprod_{s>0}
 \{(x_1,\dotsc,x_s)\mid
   x_1,\dotsc,x_s\in X_n^\x,\ x_1<\dotsc<x_s\}.
 $$
 For $i=(x_1,\dotsc,x_s)\in I$, put
 $|i|=s$ and
 $e_i=[x_1\dotso x_s]\in C_n(X^{\smash s})=M_n(X)(\<s\>)$.
 The elements $e_i$ form a basis of $M_n(X)$ in the following
 sense.
 For any
 functor $F\:\Omega^\o\to\Mod R$ and
 elements $a_i\in F(|i|)$, $i\in I$,
 there exists a unique functor morphism $T\:M_n(X)\to F$
 such that $\^{|i|}T(e_i)=a_i$ for all $i\in I$.
 Therefore,
 for a functor epimorphism $\~F\to F$,
 any functor morphism $M_n(X)\to F$ lifts to $\~F$.
 \qed


 \head {2. Main constructions}


 \subhead {2.A. Diagonal complexes.}
 A {\it bicomplex\/} $\_W^*_*$ (of $R$-modules) has
 differentials $d'\:\_W^{p-1}_q\to\_W^p_q$ and
 $d''\:\_W^p_q\to\_W^p_{q-1}$,
 which commute: $d''d'=d'd''$.
 The {\it diagonal\/} (or {\it complete total\/})
 chain complex $\diag_*\_W^*_*=W_*$ of $\_W^*_*$ has
 $$
 W_n=\prod_{q-p=n}\_W^p_q.
 $$
 For $w\in W_n$,
 we have $w=(w^p_q)_{q-p=n}$,
 where $w^p_q\in\_W^p_q$.
 The differential $\d\:W_n\to W_{n-1}$ is defined by
 $$
 (\d w)^p_q=d''(w^p_{q+1})-(-1)^nd'(w^{p-1}_q),
 \quad q-p=n-1.
 $$


 \subhead {2.B. The complex $\Hom_*(U_*,V_*)$.}
 Given chain complexes $U_*$ and $V_*$ in some $R$-linear
 category,
 we define the bicomplex $\_\Hom^*_*(U_*,V_*)$ with
 $\_\Hom^p_q(U_*,V_*)=\Hom(U_p,V_q)$ and
 the differentials induced by those of $U_*$ and $V_*$.
 We have
 $$
 \Hom_*(U_*,V_*)=\diag_*\_\Hom^*_*(U_*,V_*).
 $$


 \subhead {2.C. The complex $G_*(X,Y)$.}
 We put
 $$
 \_G^*_*(X,Y)=\_\Hom^*_*(M_*(X),M_*(Y)), \quad
 G_*(X,Y)=\Hom_*(M_*(X),M_*(Y)).
 $$


 \subhead {2.D. Construction of $\lambda_*(X,Y)$.}
 For $s>0$,
 let $\^s\eta\:Y^X\smash X^{\smash s}\to Y^{\smash s}$ be the
 evaluation map.
 For $s>0$ and $p,q\in\Z$,
 we
 have the homomorphism $C_q(\^s\eta)\:
 C_q(Y^X\smash X^{\smash s})\to C_q(Y^{\smash s})$ and
 define the homomorphism
 $$
 \^s\lambda^p_q\:C_{q-p}(Y^X)\to
 \Hom(C_p(X^{\smash s}),C_q(Y^{\smash s}))
 $$
 by
 $$
 \^s\lambda^p_q(z)(u)=C_q(\^s\eta)(z\times u),
 \quad u\in C_p(X^{\smash s}),
 \quad z\in C_{q-p}(Y^X).
 $$
 The homomorphisms $\^s\lambda^p_q$ form
 the promised chain homomorphism $\lambda_*(X,Y)$.


 \subhead {2.E. The complex $D_*(V)$.}
 For a cosimplicial space $V$,
 we have the bicomplex $\_D^*_*(V)$ with
 $\_D^p_q(V)=C_q(V^p)$ and
 the following differentials.
 The differential $d'\:C_q(V^{p-1})\to C_q(V^p)$ is defined by
 $$
 d'=\sum_{i=0}^p(-1)^iC_q(\delta^i),
 $$
 where $\delta^i\:V^{p-1}\to V^p$ are the coface maps.
 The differential $d''\:C_q(Y^{X_p})\to C_{q-1}(Y^{X_p})$ is
 the ordinary boundary operator.
 We put $D_*(V)=\diag_*\_D^*_*(V)$.


 \subhead {2.F. The complex $D_*(X,Y)$.}
 Consider the cosimplicial space $V=\hom(X,Y)$ with
 $V^p=Y^{X_p}$ \cite[Ch.\ X, 2.2 (ii)]{BouKan}.
 We put
 $$
 \_D^*_*(X,Y)=\_D^*_*(V), \quad
 D_*(X,Y)=D_*(V).
 $$


 \subhead {2.G. Construction of $\mu_*(X,Y)$.}
 For $x\in X_p$,
 we have the composite map
 $$
 \xymatrix {
 \theta^x\:
 Y^X\smash\Delta^p_+
 \ar[r]^-{\id\smash\=x} &
 Y^X\smash X
 \ar[r]^-{\eta} &
 Y,
 }
 $$
 where
 $\=x\:\Delta^p_+\to X$ is the characteristic map of the
 simplex $x$ and
 $\eta$ is the evaluation map.
 Combining $\theta^x$ over all $x\in X_p$,
 we get a map
 $$
 \theta^p\:Y^X\smash\Delta^p_+\to Y^{X_p}.
 $$
 For $p\ge0$ and $q\in\Z$,
 we
 have the homomorphism
 $C_q(\theta^p)\:C_q(Y^X\smash\Delta^p_+)\to C_q(Y^{X_p})$ and
 introduce the homomorphism
 $$
 \mu^p_q\:C_{q-p}(Y^X)\to C_q(Y^{X_p}),\quad
 \mu^p_q(z)=C_q(\theta^p)(z\times[\iota_p]).
 $$
 The homomorphisms $\mu^p_q$ form
 the promised chain homomorphism $\mu_*(X,Y)$.


 \subhead {2.H. Construction of $\epsilon_*(X,Y)$.}
 A simplex $v\in(Y^{X_p})_q$ is a basepoint-preserving function
 $v\:X_p\to Y_q$.
 For $s>0$ and $p,q\ge0$,
 we define the homomorphism
 $$
 \^s\epsilon^p_q\:C_q(Y^{X_p})\to
 \Hom(C_p(X^{\smash s}),C_q(Y^{\smash s}))
 $$
 by
 $$
 \^s\epsilon^p_q([v])([x_1\dotso x_s])=
 [v(x_1)\dotso v(x_s)],
 \quad x_1,\dotsc,x_s\in X_p,
 \quad v\in(Y^{X_p})_q.
 $$
 The homomorphisms $\^s\epsilon^p_q$ form
 a homomorphism of bicomplexes
 $$
 \_\epsilon^*_*(X,Y)\:\_D^*_*(X,Y)\to\_G^*_*(X,Y)
 $$
 and thus the promised chain homomorphism $\epsilon_*(X,Y)$.

 \begin {remark} [Remark.]
 The bicomplexes $\_D^*_*(X,Y)$ and $\_G^*_*(X,Y)$ are in fact
 cosimplicial simplicial $R$-modules.
 (To see this, recall that,
 for every space $Z$,
 $C_*(Z)$ is in fact a simplicial $R$-module and thus
 $M_*(Z)$ is a simplicial functor.)
 The homomorphism $\_\epsilon^*_*(X,Y)$ preserves this
 structure.
 \end {remark}

 One easily verifies that
 $\epsilon_*(X,Y)\circ\mu_*(X,Y)=\lambda_*(X,Y)$.


 \subhead {2.I. Proof of Theorem 0.6.}
 Take $p,q\ge0$.
 It suffices to prove that the homomorphism
 $$
 \epsilon^p_q=(\^s\epsilon^p_q)_{s>0}\:C_q(Y^{X_p})\to
 \Hom(M_p(X),M_q(Y))
 $$
 is an isomorphism.
 We
 construct a homomorphism
 $$
 \xi^p_q\:\Hom(M_p(X),M_q(Y))\to C_q(Y^{X_p})
 $$
 and
 leave to the reader to verify that $\xi^p_q\circ\epsilon^p_q$
 and $\epsilon^p_q\circ\xi^p_q$ are the identities.

 Fix a linear order on $X_p^\x$.
 Suppose we are given sets $E,F\subseteq X_p^\x$ such that
 $E\supseteq F\neq\varnothing$.
 We have $E=\{x_1,\dotsc,x_s\}$
 for some $x_1<\dotso<x_s$.
 Put $\kappa_E=x_1\dotso x_s\in(X^{\smash s})_p$.
 For $y_1,\dotsc,y_s\in Y_q$,
 define
 the function $\phi_E^F(y_1,\dotsc,y_s)\:X_p\to Y_q$
 by the rules
 \begin{itemize}
 \item[]
 $x_t\mapsto y_t$ for $t=1,\dotsc,s$ such that $x_t\in F$;
 \item[]
 $x\mapsto*$ for all other $x\in X_p$.
 \end{itemize}
 We have the homomorphism
 $\Phi_E^F\:C_q(Y^{\smash s})\to C_q(Y^{X_p})$
 with
 $\Phi_E^F([y_1\dotso y_s])=[\phi_E^F(y_1,\dotsc,y_s)]$
 for $y_1,\dotsc,y_s\in Y_q^\x$.
 Define the homomorphism
 $$
 \psi_E^F\:
 \Hom_{\Sigma_s}(C_p(X^{\smash s}),C_q(Y^{\smash s}))\to
 C_q(Y^{X_p})
 $$
 by $\psi_E^F(t)=\Phi_E^F(t([\kappa_E]))$.
 (One may note that
 $\psi_E^F$ does not depend on the order on $X_p^\x$.)
 For a functor morphism $T\:M_p(X)\to M_q(Y)$,
 we set
 $$
 \xi^p_q(T)=
 \sum_{E,F\subseteq X_p^\x:E\supseteq F\neq\varnothing}
 (-1)^{|E|-|F|}
 \psi_E^F(\^{|E|}T).
 $$
 \qed


 \head {3. Anderson's model}


 \subhead {3.A. General cosimplicial case.}
 We follow \cite[\S~2]{Bou}.
 Let $V$ be a cosimplicial space.
 We have the (unbounded) chain complex $D_*(V)$ (see 2.E).
 There is the chain homomorphism
 $$
 \mu_*(V)\:C_*(\Tot U)\to D_*(V)
 $$
 formed by the homomorphisms
 $$
 \mu^p_q\:C_{q-p}(\Tot U)\to C_q(V^p)
 $$
 that are defined in the following way.
 A simplex $w\in(\Tot V)_n$ is a sequence $(w^p)_{p\ge0}$ of
 maps $w^p\:\Delta^n_+\smash\Delta^p_+\to V^p$.
 For $w\in(\Tot U)_{q-p}$,
 we
 have the homomorphism
 $C_q(w^p)\:C_q(\Delta^{q-p}_+\smash\Delta^p_+)\to C_q(V^p)$
 and set
 $$
 \mu^p_q([w])=C_q(w^p)([\iota_{q-p}]\times[\iota_p]).
 $$
 
 \begin{claim} [3.1. Theorem.]
 Suppose that
 $R=\Z/\ell$, $\ell$ a prime,
 $V$ is fibrant and
 the spaces $V^p$, $p\ge0$, and $\Tot V$ are $\ell$-toy.
 Then $\mu_*(V)$ is a quasi-isomorphism.
 \end{claim}

 \begin{demo} [Proof.]
 Apply
 Shipley's strong convergence theorem \cite[Theorem 6.1]{Shi}
 and \cite[Lemma 2.3]{Bou}.
 \qed
 \end{demo}


 \subhead {3.B. Proof of Theorem 0.5.}
 We have
 the cosimplicial space $V=\hom(X,Y)$ and
 the canonical isomorphism $Y^X=\Tot V$
 \cite[Ch.\ X, 3.3 (i)]{BouKan}.
 The diagram
 $$
 \xymatrix {
 C_*(Y^X)
 \ar[rr]^-{\mu_*(X,Y)}
 \ar@{=}[d] & &
 D_*(X,Y)
 \ar@{=}[d] \\
 C_*(\Tot V)
 \ar[rr]^-{\mu_*(V)} & &
 D_*(V)
 }
 $$
 is commutative.

 The cosimplicial space $V$ is fibrant
 by \cite[Ch.\ X, 4.7 (ii)]{BouKan}.
 The spaces $V^p$ are $\ell$-toy since
 $X$ is gradual and
 $Y$ is $\ell$-toy.
 The spaces $Y^X$ and thus $\Tot V$ are $\ell$-toy since
 $X$ is compact and
 $Y$ is fibrant and $\ell$-toy.
 By Theorem 3.1,
 $\mu_*(V)$ is a quasi-isomorphism.
 \qed


 \head {4. Arone's model}


 \subhead {4.A. Homotopy invariance.}
 \begin {claim} [4.1. Lemma.]
 Let
 $e\:X'\to X$ and $f\:Y\to Y'$ be weak equivalences of spaces.
 Suppose that $Y$ and $Y'$ are fibrant.
 Then
 $\lambda_*(X,Y)$ is a quasi-isomorphism if and only if
 $\lambda_*(X',Y')$ is.
 \end {claim}

 \begin {demo} [Proof.]
 The maps $e$ and $f$ induce a map
 $g\:Y^X\to Y'^{\,X'}$.
 We have the commutative diagram
 $$
 \xymatrix {
 C_*(Y^X)
 \ar[rr]^-{\lambda_*(X,Y)}
 \ar[d]_{C_*(g)} & &
 G_*(X,Y)
 \ar[d]^-{G_*(e,f)} \\
 C_*(Y'^{\,X'})
 \ar[rr]^-{\lambda_*(X',Y')} & &
 G_*(X',Y').
 }
 $$
 $C_*(g)$ is a quasi-isomorphism since
 $g$ is a weak equivalence.
 It follows from Lemma 0.2 that
 $G_*(e,f)$ is a quasi-isomorphism.
 The desired equivalence is clear now.
 \qed
 \end {demo}


 \subhead {4.B. Proof of Theorem 0.3.}
 If $X$ is compact,
 the assertion follows immediately from Theorems 0.5 and 0.6.
 In general,
 $X$ is weakly equivalent to a compact space $X^\circ$.
 Using Lemma 4.1,
 we pass from $\lambda_*(X^\circ,Y)$ to $\lambda_*(X,Y)$.
 \qed


 \head {5. Reconstructing $X$ from $M_*(X)$}


 \subhead {5.A. Composition of maps and homomorphisms.}

 \begin{claim} [5.1. Lemma.]
 Let
 $X$, $Y$ and $Z$ be spaces and
 $\gamma\:Z^Y\smash Y^X\to Z^X$ be the composition map.
 Then the diagram of chain complexes and chain homomorphisms
 $$
 \xymatrix {
 C_*(Z^Y)\otimes C_*(Y^X)
 \ar[rr]^-{\text{\rm cross product}}
 \ar[d]^-{\lambda_*(Y,Z)\otimes\lambda_*(X,Y)} & &
 C_*(Z^Y\smash Y^X)
 \ar[rr]^-{C_*(\gamma)} & &
 C_*(Z^X)
 \ar[d]_-{\lambda_*(X,Z)} \\
 G_*(Y,Z)\otimes G_*(X,Y)
 \ar[rrrr]^-{\text{\rm composition}} & & & &
 G_*(X,Z)
 }
 $$
 is commutative.
 \end{claim}

 \begin{demo} []
 This follows from the associativity of the cross product.
 \qed
 \end{demo}


 \subhead {5.B. Proof of Corollary 0.4.}
 Lemma 0.2 allows us to assume $X$ and $Y$ fibrant.
 Note that
 $H_0(G_*(X,Y))=[M_*(X),M_*(Y)]$,
 the $R$-module of chain homotopy classes.
 By Lemma 5.1,
 we have the commutative diagram
 $$
 \xymatrix {
 H_0(X^Y)\otimes H_0(Y^X)
 \ar[r]^-{\text{\rm cross product}}
 \ar[d]^-{H_0(\lambda_*(Y,X))\otimes H_0(\lambda_*(X,Y))} &
 H_0(X^Y\smash Y^X)
 \ar[r]^-{H_0(\gamma)} &
 H_0(X^X)
 \ar[d]_-{H_0(\lambda_*(X,X))} \\
 [M_*(Y),M_*(X)]\otimes[M_*(X),M_*(Y)]
 \ar[rr]^-{\text{\rm composition}} & &
 [M_*(X),M_*(X)],
 }
 $$
 where $\gamma\:X^Y\smash Y^X\to X^X$ is the composition map.
 We use the notation $B\otimes A\mapsto B\circ A$ for the upper
 line homomorphism $H_0(X^Y)\otimes H_0(Y^X)\to H_0(X^X)$.
 By Theorem 0.3,
 $H_0(\lambda_*(X,Y))$,
 $H_0(\lambda_*(Y,X))$ and
 $H_0(\lambda_*(X,X))$
 are isomorphisms.

 Let $f\:M_*(X)\to M_*(Y)$ and $g\:M_*(Y)\to M_*(X)$ be
 mutually inverse chain homotopy equivalences.
 We have
 $[f]=H_0(\lambda_*(X,Y))(A)$ for some $A\in H_0(Y^X)$ and
 $[g]=H_0(\lambda_*(Y,X))(B)$ for some $B\in H_0(X^Y)$.
 By the diagram,
 $B\circ A=1$ in $H_0(X^X)$.
 Thus there are maps $a\:X\to Y$ and $b\:Y\to X$ such that
 $b\circ a\sim\id_X$.
 Interchanging $X$ and $Y$ in this reasoning,
 we get maps $a'\:X\to Y$ and $b'\:Y\to X$ such that
 $a'\circ b'\sim\id_Y$.
 Since $X$ and $Y$ are $\ell$-toy,
 these four maps are weak equivalences.
 \qed




 {\noindent \tt ssp@pdmi.ras.ru}

 {\noindent \tt http://www.pdmi.ras.ru/\tildeaccent{}ssp}

 \end {document}